\documentclass[12pt,a4paper]{article}

\newcommand{\eqref}[1]{\textup{(\ref{#1})}}
\newcommand{\text}[1]{\mbox{#1}}
\newcommand{\operatorname}[1]{\mathop{\mathrm{#1}}\nolimits}
\newcommand{\DeclareMathOperator}[2]{\newcommand{#1}{\operatorname{#2}}}
\newcommand{\binom}[2]{{#1\choose#2}}

\newcommand{\href}[2]{#2}

\IfFileExists{umsb.fd}{%
  \DeclareMathAlphabet{\mathbb}{U}{msb}{m}{n}%
}{%
  \newcommand{\mathbb}[1]{\mathbf{#1}}%
}
\IfFileExists{ueuf.fd}{%
  \DeclareMathAlphabet{\mathfrak}{U}{euf}{m}{n}%
}{%
  \newcommand{\mathfrak}[1]{\mathit{#1}}%
}

\providecommand{\C}{}
\renewcommand{\C}{\mathbb C}
\newcommand{\R}{\mathbb R}
\renewcommand{\Im}{\operatorname{Im}}
\renewcommand{\Re}{\operatorname{Re}}
\DeclareMathOperator{\Herm}{Herm}
\DeclareMathOperator{\Mat}{Mat}

\DeclareMathOperator{\rank}{rank}
\DeclareMathOperator{\im}{im}
\DeclareMathOperator{\Span}{span}
\DeclareMathOperator{\tr}{tr}

\newtheorem{theorem}{Theorem}
\newtheorem{lemma}[theorem]{Lemma}
\newtheorem{proposition}[theorem]{Proposition}

\newcommand{\proofsquare}{%
  \begingroup\setlength{\fboxsep}{0pt}%
  \fbox{\rule{0pt}{0.6em}\rule{0.6em}{0pt}}\endgroup
}
\newenvironment{proof}[1][Proof]{%
  \par\addvspace{\medskipamount}\normalfont\noindent
  \textit{#1.}\ \ignorespaces
}{%
  \unskip\nobreak\hfill\proofsquare
  \par\addvspace{\medskipamount}\ignorespacesafterend
}

\title{Description of the Lattice of Exceptional CR Types}
\author{R.~M.~Aminev\thanks{Email: \texttt{ramazan.aminev@math.msu.ru}},
	V.~K.~Beloshapka\thanks{Email: \texttt{vkb@strogino.ru}}\\
	\normalsize M.~V.~Lomonosov Moscow State University}
\date{}
\begin{document}
\maketitle
\begin{abstract}
We prove that an exceptional CR quadric of type $(n,k)$ exists
if and only if $4\le n$ and $4\le k\le(n-2)^2+1$.
Explicit examples of exceptional quadrics are given for all these types.
The proof of nonexistence for the remaining types is based on
\cite{VB22} and \cite{Tu}. This paper completes the investigation
undertaken in~\cite{VB26}.
\end{abstract}

\vspace{3ex}

\section{Introduction}

\vspace{3ex}

Let $M$ be a germ of a generic real-analytic submanifold of
$\C^{n+k}$ at a point $p$. By a translation in $\C^{n+k}$,
we may assume that $p=0$. Its CR dimension is the complex dimension
$n=\dim_\C(T_0M\cap iT_0M)$, and $k$ is the real codimension
of $M$ in $\C^{n+k}$.
We call the pair~$(n,k)$ its $CR$ type.
In suitable holomorphic coordinates, the equation of $M$ takes the form
\[
 \Im w=F(z,z)+O(3),\qquad
 z\in\C^n,\quad w\in\C^k,
\]
where $O(3)$ denotes terms of weight at least three.
Weights are determined by assigning the coordinate weights $[z]=1$, $[w]=2$.
The components of the form $F$ are represented by Hermitian matrices
$H_1,\ldots,H_k$ of size $n\times n$, that is,
$ F(z,b)=(b^*H_1z,\ldots,b^*H_kz)^T$.
The corresponding \emph{tangent quadric} is
\[
 Q=\{(z,w)\in\C^n\times\C^k:\Im w=F(z,z)\}.
\]
Here $H^*=\bar H^T$, while $H^T$ denotes transposition without
conjugation; the vectors $z,b,w$ are regarded as columns.
We denote the space of Hermitian matrices by $\mathrm{Herm}(n)$.
Let $\mathrm{aut}(M)$ be the algebra of germs of holomorphic vector fields
whose real parts are tangent to $M$.
The finite-dimensionality criterion for $\mathrm{aut}(Q)$ consists
of the following two conditions (see~\cite{VB88}):
\[
 H_1,\ldots,H_k\text{ are linearly independent over }\R,
 \qquad \bigcap_{j=1}^k\ker H_j=\{0\}.
\]
We call such a quadric \emph{nondegenerate}.
Throughout the paper, all quadrics are assumed to be nondegenerate.
Set $ \mathcal H=\Span_\R(H_1,\ldots,H_k),\; S=\Span_\C(H_1,\ldots,H_k)$.
Then
$
 \dim_\C S=\dim_\R\mathcal H=k\le n^2,\; S^*=S,
 \; S\cap\Herm(n)=\mathcal H.
$

The tangency condition for the real vector field $2\Re X$,
where $X=f\partial_z+g\partial_w$ is holomorphic, takes the form
\begin{equation}\label{eq:tangency}
 \Im g(z,w)=2\Re F(f(z,w),z),\qquad w=u+iF(z,z).
\end{equation}

If we extend the definition of weights by assigning to differentiation
with respect to $z$ and $w$ the weights $(-1)$ and $(-2)$, respectively,
then, assuming that $Q$ is nondegenerate, we obtain a grading of the
real algebra $g(Q)=\mathrm{aut}(Q)$:
\[
 g(Q)=\bigoplus_{j\ge-2}^{l} g_j(Q),\qquad
 [g_i,g_j]\subseteq g_{i+j}.
\]
We call the parameter $l\geq0$ the \emph{length}.
The first example of a quadric of length greater than two appeared
in~\cite{M} (in this example $l=6$, see~\cite{GM}).
We call quadrics whose algebras have length greater than two
\emph{exceptional}. A series of exceptional quadrics was constructed
in~\cite{GM}. The papers \cite{VB22} and \cite{Tu} showed that
exceptional quadrics do not exist for many $CR$ types.
The paper~\cite{VB26} is devoted to partitioning the lattice of $CR$ types
into exceptional types (for which exceptional quadrics exist)
and nonexceptional types (for which they do not).

The negative part $g_{-}=g_{-2}+g_{-1}$ is fundamental,
and its prolongation is transitive in Tanaka's sense
(see~\cite{Ta} and also \cite{Tu}, \S\,2).
In particular, for $j\ge0$, the condition $[X,g_{-1}]=0$, where $X\in g_j$,
implies $X=0$. Hence
\begin{equation}\label{eq:transitivity}
\mbox{if } g_3=0,\quad\mbox{then}\quad g_j=0\quad(j\ge3).
\end{equation}
Indeed, if $g_{j-1}=0$, then for $X\in g_j$ we have
$[X,g_{-1}]\subseteq g_{j-1}=0$, and therefore $X=0$.
Thus, to rule out all higher weights, it suffices to rule out weight three.
We shall prove the following theorem.

\begin{theorem}\label{th:main}
An exceptional nondegenerate CR quadric of type $(n,k)$ exists
if and only if
\begin{equation}\label{eq:iff}
 4 \le n,\qquad 4\le k\le(n-2)^2+1.
\end{equation}
\end{theorem}

For every type satisfying~\eqref{eq:iff}, there is an example
with $g_3\ne0$ and $g_4\ne0$.
For every nondegenerate quadric outside this range,
$g_j=0$ for all $j\ge3$.
Termination at weight two means that all higher weights are absent;
the component $g_2$ itself may also be zero.

\vspace{3ex}

Interest in quadrics stems primarily from the properties of their
holomorphic symmetries~\cite{Ta}, \cite{CM}, \cite{VB02}.
Their groups of holomorphic automorphisms have maximal dimension
within the class of manifolds of the simplest Bloom--Graham type~\cite{BG}
for which they serve as model manifolds:
$$\dim_\R\mathrm{aut}(M)
 \le \dim_\R\mathrm{aut}(Q).$$
For germs of arbitrary Bloom--Graham type, this role is played
by a broader class of model manifolds~\cite{VB20}.

We also note that the main technical difficulties in the proof
of the theorem arise in treating the boundary parabola
$k=(n-2)^2+2$.

\vspace{3ex}

\section{Construction of examples of exceptional quadrics}

\vspace{3ex}

\begin{proposition}\label{prop:allk}
For every $4\le n$ and $4\le k\le(n-2)^2+1$, there exists
a nondegenerate quadric of type $(n,k)$ for which $g_3\ne0$ and $g_4\ne0$.
\end{proposition}
\begin{proof}
Set \(m=n-2\) and split the coordinates into
\[
 x=(x_1,\ldots,x_m)\in\C^m,\qquad y=(y_1,y_2)\in\C^2.
\]
Choose a basis of \(\Herm(m)\) beginning with
\[
 h_1=x_1\bar x_2+x_2\bar x_1,\quad
 h_2=ix_1\bar x_2-ix_2\bar x_1,\quad h_3=|x_2|^2.
\]
Such a basis exists because these three forms are independent.
Take its first $k-1\le m^2$ elements and define the quadric
\begin{equation}\label{eq:ex-example}
\begin{array}{@{}r@{}l@{}}
\displaystyle \Im w_j & \displaystyle {}=h_j(x,\bar x),\quad 1\le j\le k-1, \\
\displaystyle \Im s & \displaystyle {}=q(x,y):=
 2\Re(x_1\bar y_1+x_2\bar y_2)+\sum_{j=3}^m|x_j|^2.
\end{array}
\end{equation}
For $m=2$, the last sum is empty.
The forms are independent because the $h_j$ are independent and $q$ is
the only form involving $y$.
The matrix of $q$ in the coordinate order $x_1,y_1,x_2,y_2,x_3,\ldots$
is
\[
 \left(\begin{array}{@{}cc@{}}0&1\\1&0\end{array}\right)
 \oplus\left(\begin{array}{@{}cc@{}}0&1\\1&0\end{array}\right)\oplus I_{m-2}.
\]
It is invertible, so the common kernel is trivial. Thus, the quadric
is nondegenerate for every chosen $k$.

\vspace{3ex}

Next, set
\[
 \Theta=2w_3x_1-(w_1-iw_2)x_2
\]
and define the vector field
\begin{equation}\label{eq:ex-X}
 X=-2iw_3\Theta\,\partial_{y_1}
    +i(w_1+iw_2)\Theta\,\partial_{y_2}.
\end{equation}
The vector field satisfies the tangency condition for the first $k-1$
equations because they do not depend on $y$.
For the last equation, we require $(X+\bar X)q=0$.

Since
$\partial_{y_1}q=\bar x_1$ and $\partial_{y_2}q=\bar x_2$, we have
\[
 Xq=\bar x_1X(y_1)+\bar x_2X(y_2)
   =i\Theta\bigl((w_1+iw_2)\bar x_2-2w_3\bar x_1\bigr).
\]
Substituting $w_j=u_j+ih_j$ and the expressions for $h_1,h_2,h_3$, we obtain
\[
\begin{array}{@{}r@{}l@{}}
\displaystyle (w_1+iw_2)\bar x_2-2w_3\bar x_1 & \displaystyle {}=(u_1+iu_2+2ix_2\bar x_1)\bar x_2
      -2(u_3+i|x_2|^2)\bar x_1 \\
\displaystyle  & \displaystyle {}=(u_1+iu_2)\bar x_2-2u_3\bar x_1, \\
\displaystyle \Theta & \displaystyle {}=2(u_3+i|x_2|^2)x_1
      -(u_1-iu_2+2ix_1\bar x_2)x_2 \\
\displaystyle  & \displaystyle {}=2u_3x_1-(u_1-iu_2)x_2.
\end{array}
\]
Thus, the expression in parentheses in $Xq$ equals $-\bar\Theta$.
Consequently, on the quadric,
\[
 Xq=-i\Theta\bar\Theta=-i|\Theta|^2,\qquad
 (X+\bar X)q=2\Re(Xq)=0.
\]
The vector field $X$ is nontrivial and has weight four: $\Theta$ has weight three,
and the coefficients of $X$ have weight five. Hence $g_4\ne0$,
and, by~\eqref{eq:transitivity}, also $g_3\ne0$.
The construction uses only $h_1,h_2,h_3$ and remains unchanged
when any chosen forms $h_4,\ldots,h_{k-1}$ are added.
\end{proof}

\vspace{3ex}

\section{Proof of the nonexistence results for exceptional quadrics}

\vspace{3ex}

On $\Mat_n(\C)$ we use the Hermitian inner product
$\tr(A^*B)$, and on $\C^n$ we use the standard
inner product $a^*b$. Orthogonal complements in both spaces are denoted by
${}^\perp$.

For a complex subspace $V\subseteq\C^n$, set
\[
 V^F=\{v\in\C^n:F(v,u)=0\text{ for all }u\in V\}.
\]
For $a\in\C^n$, we write $a^F=(\C a)^F$.
The identity $a^*H_jv=(H_ja)^*v$ implies
\begin{equation}\label{eq:annihilator}
 a^F=(Sa)^\perp.
\end{equation}

A vector $a$ is called \emph{cyclic} for $S$ if
$Sa=\C^n$, that is, the columns $H_1a,\ldots,H_ka$ span $\C^n$.

We use the following sufficient conditions for $g_3=0$, established
in Tumanov's paper~\cite{Tu} (Theorem 5.2 and Propositions 6.2 and 6.4):
\begin{eqnarray}
\exists a:\ Sa=\C^n & \Longrightarrow & g_3=0;
 \label{eq:cyclic} \\
(a^F)^F\cap a^F=0\text{ for generic }a
        & \Longrightarrow & g_3=0.\label{eq:rad}
\end{eqnarray}

Here and below, ``for generic $a$'' means on an open dense subset
of $\C^n\simeq\R^{2n}$.
We also note that the first condition is a special case of the second.

\vspace{3ex}
The absence of a cyclic vector can be interpreted in terms of linear spaces of matrices of bounded rank. Indeed, in this case, the linear space of matrices $(\mathcal M(a)=(H_1a,\ldots,H_ka), a\in\C^n)$ consists of matrices of rank less than $n$. The theory of such spaces goes back to the work of Flanders~\cite{F}, Atkinson and Lloyd~\cite{AL}; recent developments are presented, in particular, in the paper of Huang and Landsberg~\cite{HL}. We do not directly use the results of these papers here, but the nature of our argument is similar.

\vspace{3ex}

Suppose that there is no cyclic vector. Let us estimate the possible dimension of $S$.
Pass from $S$ to its orthogonal complement. Since $S^*=S$, it has the form
\begin{equation}\label{eq:trace-ann}
 N=S^\perp=\{M\in\Mat_n(\C):\tr(HM)=0\text{ for all }H\in S\},
 \qquad d=n^2-k.
\end{equation}
The nondegeneracy of the inner product gives $\dim_\C N=d$.
The equality $S^*=S$ implies $N^*=N$, since, for $M\in N$ and $H\in S$,
\[
 \tr(HM^*)=\overline{\tr(MH^*)}
          =\overline{\tr(H^*M)}=0.
\]
Hereafter, all dimensions of spaces of matrices and polynomials are taken
over $\C$, unless otherwise specified.

Thus, assume that \(Sa\ne\C^n\) for all \(a\). Consider
\[
 \mathcal M(a)=(H_1a\ \cdots\ H_ka).
\]
All its minors of order $n$ vanish identically.
Therefore, over the field $\C(a_1,\ldots,a_n)$, the matrix has rank less than $n$,
and its left kernel is nonzero. Hence there exists a polynomial
column \(p(a)\ne0\) such that
\begin{equation}\label{eq:poly-p}
 p(a)^THa=0\quad\text{for all }H\in S.
\end{equation}

By the cyclicity of the trace, $p(a)^THa=\tr\bigl(Hap(a)^T\bigr)$. It then follows from~\eqref{eq:trace-ann} and~\eqref{eq:poly-p} that
\begin{equation}\label{eq:rankone-N}
 ap(a)^T\in N.
\end{equation}
Since $Ha$ is linear in $a$, each homogeneous component of $p$
also satisfies~\eqref{eq:poly-p}. Among all nonzero homogeneous
polynomial columns satisfying~\eqref{eq:poly-p}, choose $p$ with minimal
\begin{equation}\label{min_r}
 r=\dim_\C\Span_\C\{p_1,\ldots,p_n\}.
\end{equation}
The quantity $r$ equals the dimension of the complex span of the values $p(a)$.
Denote this span by
$R=\Span_\C\{p(a):a\in\C^n\}\subseteq\C^n$.
If $r=1$, then $p(a)=f(a)v$ for a constant $v\ne0$, and
\eqref{eq:poly-p} gives $v^TH=0$ for all $H\in S$.
In particular, $v^TH_j=0$, and the Hermitian property of $H_j$ implies
$H_j\bar v=0$ for all $j$. This contradicts
the triviality of the common kernel. Thus,
\[
 2\le r\le n.
\]
Choose a unitary matrix $U$ taking $R$ to the coordinate
subspace $\C^r\times\{0\}$. Under the change of variables $a=U^Tz$,
the corresponding transformations are
\[
 H'=\bar UHU^T,\qquad p'(z)=Up(U^Tz),\qquad N'=\bar UNU^T.
\]
They preserve~\eqref{eq:poly-p} and~\eqref{eq:rankone-N}:
\[
 p'(z)^TH'z=p(U^Tz)^THU^Tz,\qquad
 zp'(z)^T=\bar U\bigl(ap(a)^T\bigr)U^T.
\]
We now return to the notation $a,p,S,N$; we have
$R=\Span_\C(e_1,\ldots,e_r)$ and $p_{r+1}=\cdots=p_n=0$.

\vspace{3ex}

Let $L$ denote the space of linear forms in $a_1,\ldots,a_n$.

\begin{lemma}\label{lem:products}
If $p_1,\ldots,p_r$ are linearly independent homogeneous polynomials
of the same degree, then
\begin{equation}\label{eq:products-bound}
 \dim(Lp_1+\cdots+Lp_r)\ge nr-\binom r2.
\end{equation}
\end{lemma}
\begin{proof}
Fix a lexicographic order on the monomials.
By Gaussian elimination, choose a basis of the span of the $p_i$ with distinct
leading monomials $m_i$. The leading monomials of its linear multiples
form the sets
\[
 A_i=\{a_1m_i,\ldots,a_nm_i\},\qquad |A_i|=n.
\]
For $i\ne j$, the intersection $A_i\cap A_j$ contains at most one
monomial: a common multiple of the next degree must coincide
with $\operatorname{\text{lcm}}(m_i,m_j)$.
Let $c(u)$ be the number of sets containing the monomial $u$.
Then
\[
 \begin{array}{@{}r@{}l@{}}
\displaystyle \left|\bigcup_i A_i\right| & \displaystyle {}=nr-\sum_u(c(u)-1)
 \ge nr-\sum_u\binom{c(u)}2 \\
\displaystyle  & \displaystyle {}=nr-\sum_{i<j}|A_i\cap A_j|
 \ge nr-\binom r2.
\end{array}
\]
Products with distinct leading monomials are linearly independent,
which gives the required estimate.
\end{proof}

\vspace{3ex}

Set
\[
 W=\Span_{\C}\{ap(a)^T:a\in\C^n\}\subset N.
\]
The linear functional $M\mapsto\sum_{i,j}c_{ij}M_{ij}$ vanishes
on $W$ if and only if $\sum_{i,j}c_{ij}a_i p_j(a)\equiv0$.
Hence the dimension of the annihilator of $W$ equals the number of independent linear
relations among the polynomials $a_i p_j(a)$, whereas $\dim W$ equals the dimension
of their complex span. By Lemma~\ref{lem:products}, we obtain
\[
 \dim W=\dim(Lp_1+\cdots+Lp_r)\ge nr-\binom r2.
\]

For matrices in $W$, only the first $r$ columns can be nonzero,
and for matrices in $W^*$, only the first $r$ rows can be nonzero.
Therefore, $W\cap W^*$ is contained in the $R\times R$ block, and
\[
 \dim(W\cap W^*)\le r^2.
\]
From \(N^*=N\), we obtain \(W+W^*\subset N\), and hence
\begin{equation}\label{eq:poly-W}
 \dim N=d\ge\dim(W+W^*)
 =2\dim W-\dim(W\cap W^*)
 \ge2\dim W-r^2.
\end{equation}

\vspace{1ex}

\begin{proposition}\label{prop:poly-bound}
Let $n\ge4$. If $S$ has no cyclic vector, then $d\ge4n-6$, that is,
\begin{equation}\label{eq:bound}
 k = \dim_\C S\le(n-2)^2+2.
\end{equation}
\end{proposition}
\begin{proof}
Consider the possible values of $r$.
\begin{enumerate}
\item If \(r=2\), then \(\dim W\ge2n-1\), whence
\[
 d\ge2(2n-1)-4=4n-6.
\]
\item If \(r=3\) and \(n\ge5\), then
\[
 d\ge2(3n-3)-9=6n-15>4n-6.
\]
\item If \(r=3,n=4\), the preceding estimate is insufficient.
In this case, \(R\) is the span of the first three coordinate vectors.
The last row of every matrix in \(W^*\) is zero, whereas the projection of \(W\)
onto the last row has dimension three:
\(a_4p_1,a_4p_2,a_4p_3\) are independent.
Therefore, the intersection \(W\cap W^*\) lies in the kernel of this projection:
\[
 \dim(W\cap W^*)\le\dim W-3.
\]
We obtain \(d\ge\dim W+3\ge9+3=12>10=4n-6\).
\item If \(r\ge4\), it suffices to take four independent components:
\[
 d\ge\dim W\ge4n-\binom42=4n-6.
\]
\end{enumerate}
In all cases, $d\ge4n-6$, which is equivalent to \eqref{eq:bound}.
\end{proof}
For $k>(n-2)^2+2$, this proposition guarantees the existence of
a cyclic vector, and condition \eqref{eq:cyclic} gives $g_3=0$.

\vspace{3ex}

It remains to consider the case $n\ge4$, $k=(n-2)^2+2$, and $d = 4n - 6$ with no cyclic vector.
In the proof of Proposition~\ref{prop:poly-bound}, the case $r=3$
gives the strict inequality $d>4n-6$ and is therefore excluded.
For $r\ge5$, applying Lemma~\ref{lem:products} to five
components gives $d\ge5n-10>4n-6$, since $n\ge r\ge5$.

Let $r=4$. Then $4n-6\le\dim W\le d=4n-6$, whence
$W=N=W^*$. Therefore, the only rows and columns of matrices in $W$
that can be nonzero are those of the $R\times R$ block, and their images lie in $R$.
On the other hand,
\begin{equation}\label{eq:images}
 \Span_\C\{\im M:M\in W\}=\C^n:
\end{equation}
when $p(a)\ne0$, the image of $ap(a)^T$ equals $\C a$, and such $a$
form a dense set. Consequently, $n=\dim R=4$.
Only the following cases remain:
\[
 r=2\qquad\text{or}\qquad(n,r,k)=(4,4,6).
\]

\vspace{3ex}

We begin with the branch $r=2$.
The inequalities $\dim W\ge2n-1$ and $4n-6=d\ge2\dim W-4$
imply
\[
 \dim(Lp_1+Lp_2)=2n-1.
\]
Since $\dim Lp_i=n$, we have $\dim(Lp_1\cap Lp_2)=1$.
Consequently, there exists a nontrivial relation
\(\ell_1p_1=\ell_2p_2\) with linear \(\ell_i\).
Divide \(p_1,p_2\) by their common factor \(\mu\):
\(p_i=\mu q_i\), \(\operatorname{\text{gcd}}(q_1,q_2)=1\).
The equality \(\ell_1q_1=\ell_2q_2\) implies \(q_1\mid\ell_2\) and
\(q_2\mid\ell_1\). Therefore, the degrees of \(q_i\) are at most one.
The independence of $p_1,p_2$ excludes degree zero.
Thus, after cancelling $\mu$, the components are linear
and independent, and identity \eqref{eq:poly-p} is preserved.
Denote these components by $\ell_1,\ell_2$ and place them
in the last two coordinates:
\[
 p(a)=(0,\ldots,0,\ell_1(a),\ell_2(a))^T.
\]
Let the last two rows of \(H\), applied to \(a\), be the
linear forms \(u_H(a),v_H(a)\).
The identity \(p(a)^THa=0\) means that
\[
 \ell_1u_H+\ell_2v_H=0.
\]
Since $\ell_1,\ell_2$ are relatively prime, it follows that
\[
 u_H=t_H\ell_2,\qquad v_H=-t_H\ell_1
\]
with a constant $t_H$ depending linearly on $H$.
Hence, the last two rows of $H$ are $t_HD$, where the rows
of the fixed matrix $D$ are the coefficients of the forms $\ell_2,-\ell_1$.
In particular, $\rank D=2$.

\vspace{3ex}

Write \(\C^n=E\oplus R\), \(\dim E=m=n-2\), \(\dim R=2\), and
\[
 D=(B\ C),\qquad C\in\Mat_2(\C).
\]
On the real vector space \(\mathcal H\), consider the map
\(H\mapsto t_H\in\C\). Its kernel consists of matrices whose
last rows, and hence, by Hermitian symmetry, last columns, are zero.
Therefore, the dimension of the kernel is at most $m^2$, and the dimension
of the image is at most two over $\R$.
Since $\dim_\R\mathcal H=m^2+2$, both bounds are attained:
the kernel is the entire upper block $\Herm(E)$, and $t_H$ ranges over $\C$.
The lower right block equals \(t_HC\) and must be Hermitian
for every \(t_H\). Substituting \(t_H=1\) and \(t_H=i\), we obtain
\[
 C=C^*,\qquad iC=(iC)^*=-iC^*.
\]
Consequently, $C=0$ and $\rank B=\rank D=2$.
Hermitian symmetry determines the upper right block, while the upper left block
is arbitrary, since the kernel contains all of $\Herm(E)$. Thus,
\begin{equation}\label{eq:block}
 \mathcal H_B=
 \left\{
 \left(\begin{array}{@{}cc@{}}A&\bar\alpha B^*\\ \alpha B&0\end{array}\right):
 A=A^*,\ \alpha\in\C
 \right\},\qquad \rank B=2.
\end{equation}
By choosing two vectors in \(E\) that \(B\) maps to the standard
basis of \(\C^2\), and adjoining a basis of \(\ker B\), we reduce
\(B\) to \((I_2\ 0)\). This is an admissible complex change of coordinates,
which acts on the forms as follows: \(H\mapsto C_0^*HC_0\).

\vspace{3ex}

The coefficients of a vector field of weight three have the form
(see~\cite{VB22}, Proposition~8;
\cite{Tu}, the proof of Theorem~5.2)
\begin{equation}\label{eq:full-weight3}
 \begin{array}{@{}r@{}l@{}}
\displaystyle f(z,w) & \displaystyle {}=A(z,z,w)+\mathcal B(w,w), \\
\displaystyle g(z,w) & \displaystyle {}=2iF\bigl(z,\mathcal B(\bar w,\bar w)\bigr),
\end{array}
\end{equation}
where $A$ is a complex trilinear map symmetric
in its first two arguments, and $\mathcal B$ is a complex symmetric
bilinear map.\\
The tangency condition \eqref{eq:tangency} yields the equation
\[
 F(A(z,z,w),b)=4iF\bigl(z,\mathcal B(\bar w,F(b,z))\bigr).
\]
Substitute $w=e_l$, where $e_l$ is a standard basis vector of $\C^k$, and set
$A_l(z)=A(z,z,e_l)$, $\mathcal B_{l\nu}=\mathcal B(e_l,e_\nu)$.
For the component with matrix $H=H_j$, we obtain
\[
 \begin{array}{@{}r@{}l@{}}
\displaystyle b^*H A_l(z) & \displaystyle {}=4i\biggl(\sum_{\nu=1}^k F_\nu(b,z)\mathcal B_{l\nu}\biggr)^*Hz \\
\displaystyle  & \displaystyle {}=\sum_{\nu=1}^k (4i\,\mathcal B_{l\nu}^*Hz)\,b^*H_{\nu}z,
\end{array}
\]
since $\overline{F_\nu(b,z)}=F_\nu(z,b)=b^*H_{\nu}z$.
Set $D_{l\nu}=4i\,\overline{\mathcal B_{l\nu}}$. Factoring out $b^*$ on the left
on both sides and using the arbitrariness of $b$, we obtain, for each $H_j$,
and by linearity for all $H\in S$,
\begin{equation}\label{eq:weight3}
 H A_l(z)=\sum_{\nu=1}^k (D_{l\nu}^T H z) H_{\nu}z,\qquad
 D_{l\nu}=D_{\nu l},\quad H\in S.
\end{equation}
If all $D_{lj}$ and $A_l$ are zero, then $\mathcal B=0$ and
$A(z,z,w)=\sum_l w_l A_l(z)=0$. By \eqref{eq:full-weight3}, the vector field is then also zero.

\vspace{3ex}

\begin{lemma}\label{lem:block}
For the quadric with the space of forms \eqref{eq:block}, we have $g_3=0$.
\end{lemma}
\begin{proof}
Let $\Phi(z,b)$ be the column vector with components $b^TH_jz$.
For $\lambda,\mu\in\C^k$, set $A_\lambda(z)=\sum_l\lambda_l A_l(z)$
and $D(\lambda,\mu)=\sum_{l,j}\lambda_l\mu_j D_{lj}$.
Multiplying \eqref{eq:weight3} on the left by $b^T$ and summing over $l$ with weights
$\lambda_l$, we obtain
\[
 b^TH A_\lambda(z)=D\bigl(\lambda,\Phi(z,b)\bigr)^THz,\qquad H\in S.
\]
Take two pairs $(z_i,b_i)$, where $z_i=(x_i,y_i)$, $b_i=(\eta_i,\theta_i)$,
$i=1,2$. Write $\varphi_i=\Phi(z_i,b_i)$ and $(u,v)=D(\varphi_1,\varphi_2)$.
By the symmetry of $D$, substituting $(z_1,b_1,\varphi_2)$ and $(z_2,b_2,\varphi_1)$
for $(z,b,\lambda)$ gives conditions on the same vector $(u,v)$.
We will show that the presence of the upper block forces $u$ to lie in both
$\C\eta_1$ and $\C\eta_2$, while the lower left block then gives
$v^TBx_1=v^TBx_2=0$. For generic pairs, this implies $(u,v)=0$.

Let $\alpha_1,\alpha_2$ denote the upper components of
$A_{\varphi_2}(z_1),A_{\varphi_1}(z_2)$, respectively.

The space $S$ contains the matrices
$H_M=\left({\setlength{\arraycolsep}{2pt}\renewcommand{\arraystretch}{0.65}\begin{array}{@{}cc@{}}\scriptstyle M & \scriptstyle 0\\\scriptstyle 0 & \scriptstyle 0\end{array}}\right)$
with arbitrary $M\in\Mat_m(\C)$ and
$H_L=\left({\setlength{\arraycolsep}{2pt}\renewcommand{\arraystretch}{0.65}\begin{array}{@{}cc@{}}\scriptstyle 0 & \scriptstyle 0\\\scriptstyle B & \scriptstyle 0\end{array}}\right)$.
For $H=H_M$ and each $i=1,2$, we have
\[
 \eta_i^TM\alpha_i=u^TMx_i\quad\text{for all }M,\qquad
 \text{hence}\quad \alpha_i\eta_i^T=x_i u^T
\]
by comparing the coefficients of the entries of $M$. For generic pairs,
the vectors $\eta_1,\eta_2$ are linearly independent, as are $Bx_1,Bx_2$,
since $m\ge2$ and $\rank B=2$. In particular, $x_i\ne0$.
Choosing a row corresponding to a nonzero coordinate of $x_i$ in the matrix equality,
we obtain $u\in\C\eta_i$ for each $i$. The intersection of these two lines
is zero. Therefore $u=0$, and since $\eta_i\ne0$, we also have $\alpha_1=\alpha_2=0$.

For $H=H_L$, the same two substitutions give
\[
 v^TBx_i=\theta_i^TB\alpha_i=0,\qquad i=1,2.
\]
The vectors $Bx_1,Bx_2$ form a basis of $\C^2$, whence $v=0$.
Thus $D(\varphi_1,\varphi_2)=0$ for generic pairs and, by polynomiality,
for all pairs. The values of $\Phi(z,b)$ span $\C^k$, since
the identity $\sum_j c_j b^TH_jz=0$ in $z,b$ implies $\sum_j c_jH_j=0$
and, by the complex linear independence of the matrices $H_j$, all $c_j=0$.
By bilinearity, $D=0$. Now \eqref{eq:weight3} implies
$H A_l(z)=0$ for all $H\in S$. The common kernel of the forms is zero, so
$A_l=0$ for each $l$, proving the lemma.
\end{proof}

This completes the argument for $r=2$. For $n\ge5$, this was
the only possibility in the boundary case.

\vspace{3ex}

\begin{lemma}\label{lem:alt4-vanishing}
For every nondegenerate quadric $Q$ of type $(4,6)$, we have $g_3(Q)=0$.
\end{lemma}
\begin{proof}
If a cyclic vector exists, the assertion follows from \eqref{eq:cyclic}.
Suppose that $Sa\ne\C^4$ for all $a$.
Consider the matrix $\mathcal M(a)=(H_1a\ \cdots\ H_6a)$ over
the field $K=\C(a_1,\ldots,a_4)$ and the homogeneous relation
$p(a)^T\mathcal M(a)=0$ with minimal $r$ chosen
above. Since $d=16-6=10$, as shown above,
only $r=2$ and $r=4$ are possible.
For $r=2$, the assertion follows from Lemma~\ref{lem:block}. Let $r=4$.

We will show that $\rank_K\mathcal M=3$. There is no cyclic vector, so the rank is at most three. Suppose that the rank is at most two;
then $\dim_K\ker\mathcal M^T\ge2$. The restriction of the functional
$v\mapsto v_4$ to this kernel has rank at most one. Hence
$\ker\mathcal M^T$ contains a nonzero vector $(q_1,q_2,q_3,0)^T$
over $K$. Multiplying by a common denominator and taking
a nonzero homogeneous component, we obtain a polynomial relation
with at most three independent components. This contradicts
the minimality of $r=4$.

Consequently, for generic $a$, we have $\dim_{\C}Sa=3$ and $p(a)\ne0$.
Set $c(a)=\overline{p(a)}$. From
$c(a)^*H_ja=p(a)^TH_ja=0$ for all $j$ and $\dim Sa=3$, we obtain
\[
 a^F=(Sa)^\perp=\C c(a).
\]

Let us verify that $F(c(a),c(a))\ne0$ for generic $a$.
Fix a nonzero matrix $H_j$ and write
$p(a)=\sum_\alpha v_\alpha a^\alpha$. The column vectors $v_\alpha$
span $\C^4$, since a nonzero linear functional vanishing on all
$v_\alpha$ would yield a linear relation between the four independent
components of $p$. If $c(a)^*H_jc(a)$ were identically zero,
then the equality
\[
 0=p(a)^TH_j\overline{p(a)}
   =\sum_{\alpha,\beta}v_\alpha^TH_j\overline{v_\beta}\,
      a^\alpha\bar a^\beta
\]
and the linear independence of the monomials $a^\alpha\bar a^\beta$ would imply
$v_\alpha^TH_j\overline{v_\beta}=0$ for all $\alpha,\beta$.
Since both $v_\alpha$ and $\overline{v_\beta}$ span $\C^4$,
this would mean that $H_j=0$, a contradiction.

For generic $a$, every vector in $(a^F)^F\cap a^F$ has the form
$\lambda c(a)$ and satisfies
\[
 0=F(\lambda c(a),c(a))=\lambda F(c(a),c(a)).
\]
Hence $\lambda=0$ and $(a^F)^F\cap a^F=\{0\}$, and condition \eqref{eq:rad} gives $g_3(Q)=0$.
\end{proof}

\vspace{3ex}

The absence of exceptional quadrics for $k\le3$ was established
in~\cite{VB22}. For types $(3,3)$ and $(3,4)$, it was also
proved there in Theorem~7. The technique developed above allows us
to exclude all cases of small CR dimension directly.

\begin{lemma}\label{lem:small-n}
For every nondegenerate quadric $Q$ of type $(n,k)$
with $n\le3$, we have $g_j(Q)=0$ for all $j\ge3$.
Moreover, if $n\le3$ and $k\ge4$, the space $S$
has a cyclic vector.
\end{lemma}
\begin{proof}
For $k\le3$, the first assertion follows from~\cite{VB22}.
Suppose that $k\ge4$. Since $k\le n^2$, the case $n=1$
is impossible. If $n=2$, we necessarily have $k=4$, so
$S=\Mat_2(\C)$ and every nonzero vector is cyclic.

It remains to consider $n=3$. Suppose that there is no cyclic vector,
and choose a homogeneous polynomial column vector $p$
satisfying~\eqref{eq:poly-p}, as above.
The triviality of the common kernel excludes $r=1$, so $r=2$ or $r=3$.
As before, set
\[
 W=\Span_\C\{ap(a)^T:a\in\C^3\}\subset N,
 \qquad d=\dim_\C N=9-k.
\]
If $r=2$, then Lemma~\ref{lem:products} gives
$\dim W\ge2\cdot3-1=5$, while
$\dim(W\cap W^*)\le4$. Since $N^*=N$, we obtain
\[
 d\ge\dim(W+W^*)
   =2\dim W-\dim(W\cap W^*)\ge10-4=6.
\]
If $r=3$, the same lemma gives
\[
 d\ge\dim W\ge3\cdot3-\binom32=6.
\]
In both cases, $k=9-d\le3$, contrary to $k\ge4$.

Thus, for $n\le3$ and $k\ge4$, a cyclic vector exists.
By~\eqref{eq:cyclic}, this implies $g_3(Q)=0$,
and by~\eqref{eq:transitivity}, we also have $g_j(Q)=0$ for all $j\ge3$.
\end{proof}

\begin{proof}[Completion of the proof of Theorem~\ref{th:main}]

Proposition~\ref{prop:allk} provides examples for all
$n\ge4$, $4\le k\le(n-2)^2+1$.
For $k\le3$, exceptionality is impossible by~\cite{VB22},
and for $n\le3$, by Lemma~\ref{lem:small-n}.
It therefore remains to consider $n\ge4$ and
$k\ge(n-2)^2+2$.
If $S$ has a cyclic vector, then $g_3=0$
by \eqref{eq:cyclic}. If it does not,
Proposition~\ref{prop:poly-bound} leaves only the equality
$k=(n-2)^2+2$. For $r=2$, the vanishing of $g_3$ follows from
Lemma~\ref{lem:block}. The remaining branch is $(n,r,k)=(4,4,6)$,
for which $g_3=0$ by Lemma~\ref{lem:alt4-vanishing}.
Thus, in all cases outside the range~\eqref{eq:iff},
we have $g_3=0$, and by~\eqref{eq:transitivity}, all components
of weight greater than two also vanish.
Together with the examples constructed above, this proves the theorem.
\end{proof}

\vspace{3ex}
The authors used ChatGPT in preparing this article.

\end{document}